# The various properties of certain subclasses of analytic functions of complex order

Nizami MUSTAFA

Department of Mathematics, Faculty of Science and Letters,
Kafkas University, Kars-TURKEY

nizamimustafa@gmail.com

31 March 2017

**Abstract**

In this paper, we introduce and investigate two new subclasses of analytic functions in the open unit disk in the complex plane. Several interesting properties of the functions belonging to these classes are examined. Here, sufficient, and necessary and sufficient conditions for the functions belonging to these classes are also given. Moreover, coefficient bound estimates for the functions belonging to these classes are obtained.

## 1. Introduction and preliminaries

Let $A$ be the class of analytic functions $f(z)$ in the open unit disk $U = \{z \in \mathbb{C} : |z| < 1\}$, normalized by $f(0) = 0 = f'(0) - 1$ of the form

$$f(z) = z + a_2 z^2 + a_3 z^3 + \cdots + a_n z^n + \cdots = z + \sum_{n=2}^{\infty} a_n z^n, a_n \in \mathbb{C}, \quad (1.1)$$

and $S$ denote the class of all functions in $A$ which are univalent in $U$.

Also, let us define by $T$ the subclass of all functions $f(z)$ in $A$ of the form

$$f(z) = z - a_2 z^2 - a_3 z^3 - \cdots - a_n z^n - \cdots = z - \sum_{n=2}^{\infty} a_n z^n, a_n \geq 0. \quad (1.2)$$

Furthermore, we will denote by $S^*(\alpha)$ and $C(\alpha)$ the subclasses of $S$ that are, respectively, starlike and convex functions of order $\alpha (\alpha \in [0,1))$. By definition (see for details, [4, 5], also [10])

$$S^*(\alpha) = \left\{ f \in S : \operatorname{Re}\left(\frac{zf'(z)}{f(z)}\right) > \alpha,\ z \in U \right\},\ \alpha \in [0,1), \quad (1.3)$$

and

$$C(\alpha) = \left\{ f \in S : \operatorname{Re}\left(1 + \frac{zf''(z)}{f'(z)}\right) > \alpha,\ z \in U \right\},\ \alpha \in [0,1). \qquad (1.4)$$

For convenience, $S^* = S^*(0)$ and $C = C(0)$ are, respectively, starlike and convex functions in $U$. It is easy to verify that $C \subset S^* \subset S$. For details on these classes, one could refer to the monograph by Goodman [5].

Note that, we will use $TS^*(\alpha) = T \cap S^*(\alpha)$, $TC(\alpha) = T \cap C(\alpha)$ and in the special case, we have $TS^* = T \cap S^*$, $TC = T \cap C$ for $\alpha = 0$.

An interesting unification of the functions classes $S^*(\alpha)$ and $C(\alpha)$ is provided by the class $S^*C(\alpha,\beta)$ of functions $f \in S$, which also satisfies the following condition:

$$\operatorname{Re}\left\{ \frac{zf'(z) + \beta z^2 f''(z)}{\beta z f'(z) + (1-\beta) f(z)} \right\} > \alpha,\ \alpha \in [0,1),\ \beta \in [0,1],\ z \in U.$$

In special case, for $\beta = 0$ and $\beta = 1$, respectively, we have $S^*C(\alpha,0) = S^*(\alpha)$ and $S^*C(\alpha,1) = C(\alpha)$, in terms of the simpler classes $S^*(\alpha)$ and $C(\alpha)$, defined by (1.3) and (1.4), respectively. Also, we will use $TS^*C(\alpha,\beta) = T \cap S^*C(\alpha,\beta)$.

The class $TS^*C(\alpha,\beta)$ was investigated by Altintaş et al. [2] and [3] (in a more general way $T_n S^* C(p,\alpha,\beta)$) and (subsequently) by Irmak et al. [6]. In particular, the class $T_n S^* C(1,\alpha,\beta)$ was considered earlier by Altintaş [1].

Inspired by the aforementioned works, we define a subclass of analytic functions as follows.

**Definition 1.1.** A function $f \in S$ given by (1.1) is said to be in the class $S^*C(\alpha,\beta;\tau)$, $\alpha \in [0,1), \beta \in [0,1], \tau \in \mathbb{C}^* = \mathbb{C} - \{0\}$ if the following condition is satisfied

$$\operatorname{Re}\left\{ 1 + \frac{1}{\tau}\left[ \frac{zf'(z) + \beta z^2 f''(z)}{\beta z f'(z) + (1-\beta) f(z)} - 1 \right] \right\} > \alpha,\ z \in U,\ \alpha \in [0,1),\ \beta \in [0,1],\ \tau \in \mathbb{C}^*.$$

In special case, we have $S^*C(\alpha,\beta;1) = S^*C(\alpha,\beta)$ for $\tau = 1$. Note that, we will use $TS^*C(\alpha,\beta;\tau) = T \cap S^*C(\alpha,\beta;\tau)$. Also, we have $TS^*C(\alpha,\beta;1) = TS^*C(\alpha,\beta)$.

In this paper, two new subclasses $S^*C(\alpha,\beta;\tau)$ and $TS^*C(\alpha,\beta;\tau)$ of the analytic functions in the open unit disk are introduced. Various characteristic properties of the functions belonging to these classes are examined. Sufficient conditions for the analytic functions belonging to the class $S^*C(\alpha,\beta;\tau)$, and necessary and sufficient conditions for those belonging to the class $TS^*C(\alpha,\beta;\tau)$ are also given. Furthermore, some coefficient bound estimates for the functions belonging to these classes are obtained.

## 2. Coefficient bounds for the classes $S^*C(\alpha,\beta;\tau)$ and $TS^*C(\alpha,\beta;\tau)$

In this section, we will examine some inclusion results of the subclasses $S^*C(\alpha,\beta;\tau)$ and $TS^*C(\alpha,\beta;\tau)$ of analytic functions in the open unit disk. Alsao, we give coefficient bound estimates for the functions belonging to these classes.

A sufficient condition for the functions in class $S^*C(\alpha,\beta;\tau)$ is given by the following theorem.

**Theorem 2.1.** Let $f \in A$. Then, the function $f(z)$ belongs to the class $S^*C(\alpha,\beta;\tau)$, $\alpha \in [0,1), \beta \in [0,1], \tau \in \mathbb{C}^* = \mathbb{C}-\{0\}$ if the following condition is satisfied

$$\sum_{n=2}^{\infty}\left[n+(1-\alpha)|\tau|-1\right]\left[1+(n-1)\beta\right]|a_n| \leq (1-\alpha)|\tau|.$$

The result is sharp for the functions

$$f_n(z) = z + \frac{(1-\alpha)|\tau|}{\left[n+(1-\alpha)|\tau|-1\right]\left[1+(n-1)\beta\right]} z^n, \; z \in U, \; n=2,3,\ldots.$$

**Proof.** According to Definition 1.1, a function $f(z)$ is in the class $S^*C(\alpha,\beta;\tau)$, $\alpha \in [0,1), \beta \in [0,1], \tau \in \mathbb{C}^* = \mathbb{C}-\{0\}$ if and only if

$$\text{Re}\left\{1+\frac{1}{\tau}\left[\frac{zf'(z)+\beta z^2 f''(z)}{\beta zf'(z)+(1-\beta)f(z)}-1\right]\right\} > \alpha.$$

It suffices to show that

$$\left|\frac{1}{\tau}\left[\frac{zf'(z)+\beta z^2 f''(z)}{\beta zf'(z)+(1-\beta)f(z)}-1\right]\right| < 1-\alpha. \tag{2.1}$$

Considering (1.1), by simple computation, we write

$$\left|\frac{1}{\tau}\left[\frac{zf'(z)+\beta z^2 f''(z)}{\beta zf'(z)+(1-\beta)f(z)}-1\right]\right| = \left|\frac{1}{\tau} \frac{\sum_{n=2}^{\infty}(n-1)[1+\beta(n-1)]a_n z^n}{z+\sum_{n=2}^{\infty}[1+\beta(n-1)]a_n z^n}\right| \leq \frac{\sum_{n=2}^{\infty}(n-1)[1+\beta(n-1)]|a_n|}{|\tau|\left\{1-\sum_{n=2}^{\infty}[1+\beta(n-1)]|a_n|\right\}}.$$

Last expression is bounded above by $1-\alpha$ if

$$\sum_{n=2}^{\infty}(n-1)[1+\beta(n-1)]|a_n| \leq |\tau|(1-\alpha)\left\{1-\sum_{n=2}^{\infty}[1+\beta(n-1)]|a_n|\right\},$$

which is equivalent to

$$\sum_{n=2}^{\infty}\left[n+(1-\alpha)|\tau|-1\right][1+(n-1)\beta]|a_n| \leq (1-\alpha)|\tau|. \tag{2.2}$$

Hence, the inequality (2.1) is true if the condition (2.2) is satisfied.

Thus, the proof of Theorem 2.1 is completed.

Setting $\tau = 1$ in Theorem 2.1, we arrive at the following corollary.

**Corollary 2.1.** The function $f(z)$ definition by (1.1) belongs to the class $S^*C(\alpha,\beta)$, $\alpha \in [0,1), \beta \in [0,1]$ if the following condition is satisfied

$$\sum_{n=2}^{\infty}(n-\alpha)[1+\beta(n-1)]|a_n| \leq 1-\alpha.$$

The result is sharp for the functions

$$f_n(z) = z + \frac{1-\alpha}{(n-\alpha)[1+\beta(n-1)]}z^n, \ z \in U, \ n = 2,3,\ldots.$$

**Remark 2.1.** The result obtained in Corollary 2.1 verifies to Corollary 2.1 in [8].

Choose $\beta = 0$ in Corollary 2.1, we have the following result.

**Corollary 2.2.** (see [9, p. 110, Theorem 1]) The function $f(z)$ definition by (1.1) belongs to the class $S^*(\alpha)$, $\alpha \in [0,1)$ if the following condition is satisfied

$$\sum_{n=2}^{\infty}(n-\alpha)|a_n| \leq 1-\alpha.$$

The result is sharp for the functions

$$f_n(z) = z + \frac{1-\alpha}{n-\alpha}z^n, \ z \in U, \ n = 2,3,\ldots.$$

Taking $\beta = 1$ in Corollary 2.1, we arrive at the following result.

**Corollary 2.3.** (see [9, p. 110, Corollary of Theorem1]) The function $f(z)$ definition by (1.1) belongs to the class $C(\alpha)$, $\alpha \in [0,1)$ if the following condition is satisfied

$$\sum_{n=2}^{\infty} n(n-\alpha)|a_n| \leq 1-\alpha.$$

The result is sharp for the functions

$$f_n(z) = z + \frac{1-\alpha}{n(n-\alpha)} z^n, \ z \in U, \ n = 2,3,\ldots.$$

**Remark 2.2.** The results obtained in Corollary 2.2 and 2.3 verifies to Corollary 2.3 and 2.4 in [7], respectively.

From the following theorem, we see that for the functions in the class $TS^*C(\alpha,\beta;\tau), \alpha \in [0,1), \beta \in [0,1], \tau \in \mathbb{R}^* = \mathbb{R}-\{0\}$ the converse of Theorem 2.1 is also true.

**Theorem 2.2.** Let $f \in T$. Then, the function $f(z)$ belongs to the class $TS^*C(\alpha,\beta;\tau)$, $\alpha \in [0,1), \beta \in [0,1], \tau \in \mathbb{R}^* = \mathbb{R}-\{0\}$ if and only if

$$\sum_{n=2}^{\infty} \left[n+(1-\alpha)|\tau|-1\right]\left[1+(n-1)\beta\right] a_n \leq (1-\alpha)|\tau|.$$

The result is sharp for the functions

$$f_n(z) = z - \frac{(1-\alpha)|\tau|}{\left[n+(1-\alpha)|\tau|-1\right]\left[1+(n-1)\beta\right]} z^n, \ z \in U, \ n = 2,3,\ldots..$$

**Proof.** The proof of the sufficiency of the theorem can be proved similarly to the proof of Theorem 2.1. Therefore, we will prove only the necessity of the theorem.

Assume that $f \in TS^*C(\alpha,\beta;\tau)$, $\alpha \in [0,1), \beta \in [0,1], \tau \in \mathbb{R}^* = \mathbb{R}-\{0\}$; that is,

$$\text{Re}\left\{1+\frac{1}{\tau}\left[\frac{zf'(z)+\beta z^2 f''(z)}{\beta z f'(z)+(1-\beta)f(z)}-1\right]\right\} > \alpha, \ z \in U. \qquad (2.3)$$

Using (1.2) and (2.3), we can easily show that

$$\text{Re}\left\{\frac{-\sum_{n=2}^{\infty}(n-1)\left[1+\beta(n-1)\right]a_n z^n}{\tau\left\{z-\sum_{n=2}^{\infty}\left[1+\beta(n-1)\right]a_n z^n\right\}}\right\} > \alpha-1.$$

The expression

$$\frac{-\sum_{n=2}^{\infty}(n-1)[1+\beta(n-1)]a_n z^n}{\tau\left\{z-\sum_{n=2}^{\infty}[1+\beta(n-1)]a_n z^n\right\}}$$

is real if choose $z$ real.

Thus, from the previous inequality letting $z \to 1$ through real values, we have

$$\frac{-\sum_{n=2}^{\infty}(n-1)[1+\beta(n-1)]a_n}{\tau\left\{1-\sum_{n=2}^{\infty}[1+\beta(n-1)]a_n\right\}} \geq \alpha - 1. \tag{2.4}$$

We will examine of the last inequality depending on the different cases of the sing of the parameter $\tau$ as follows.

Let us $\tau > 0$. Then, from (2.4), we have

$$-\sum_{n=2}^{\infty}(n-1)[1+\beta(n-1)]a_n \geq (\alpha-1)\tau\left\{1-\sum_{n=2}^{\infty}[1+\beta(n-1)]a_n\right\},$$

which is equivalent to

$$\sum_{n=2}^{\infty}[n+(1-\alpha)\tau-1][1+\beta(n-1)]a_n \leq (1-\alpha)\tau. \tag{2.5}$$

Now, let us $\tau < 0$. Then, since $\tau = -|\tau|$, from (2.4), we get

$$\frac{\sum_{n=2}^{\infty}(n-1)[1+\beta(n-1)]a_n}{|\tau|\left\{1-\sum_{n=2}^{\infty}[1+\beta(n-1)]a_n\right\}} \geq \alpha - 1;$$

that is,

$$\sum_{n=2}^{\infty}(n-1)[1+\beta(n-1)]a_n \geq (\alpha-1)|\tau|\left\{1-\sum_{n=2}^{\infty}[1+\beta(n-1)]a_n\right\}.$$

Therefore,

$$\sum_{n=2}^{\infty}[n+(\alpha-1)|\tau|-1][1+\beta(n-1)]a_n \geq -(1-\alpha)|\tau|.$$

Since $\alpha < 1$ (or $1-\alpha > \alpha-1$), from the last inequality, we have

$$\sum_{n=2}^{\infty}\left[n+(1-\alpha)|\tau|-1\right]\left[1+\beta(n-1)\right]a_n \geq -(1-\alpha)|\tau|. \tag{2.6}$$

Thus, from (2.6) and (2.7), the proof of the necessity of theorem; that is, the proof of theorem is completed.

Special case of Theorem 2.2 has been proved by Altintaş *et al* [2], $\tau=1$ (there $p=n=1$).

Setting $\tau=1$ in Theorem 2.2, we arrive at the following corollary.

**Corollary 2.4.** The function $f(z)$ definition by (1.2) belongs to the class $TS^*C(\alpha,\beta)$, $\alpha \in [0,1), \beta \in [0,1]$ if and only if

$$\sum_{n=2}^{\infty}(n-\alpha)\left[1+\beta(n-1)\right]a_n \leq 1-\alpha.$$

**Remark 2.3.** The result obtained in Corollary 2.4 verifies to Theorem 1 in [2].

Taking $\beta=0$ in Corollary 2.4, we have the following result.

**Corollary 2.5.** (see [9, p. 110, Theorem 2]) The function $f(z)$ definition by (1.2) belongs to the class $TS^*(\alpha)$, $\alpha \in [0,1)$ if and only if

$$\sum_{n=2}^{\infty}(n-\alpha)a_n \leq 1-\alpha.$$

Choose $\beta=1$ in Corollary 2.4, we arrive at the following result.

**Corollary 2.6.** (see [9, p. 111, Corollary 2]) The function $f(z)$ definition by (1.2) belongs to the class $TC(\alpha)$, $\alpha \in [0,1)$ if and only if

$$\sum_{n=2}^{\infty}n(n-\alpha)a_n \leq 1-\alpha.$$

**Remark 2.4.** The results obtained in Corollary 2.5 and 2.6 verifies to Corollary 2.7 and 2.8 in [7], respectively.

On the coefficient bound estimates of the functions belonging in the class $TS^*C(\alpha,\beta;\tau)$, we give the following result.

**Theorem 2.3.** Let the function $f(z)$ definition by (1.2) belongs to the class $TS^*C(\alpha,\beta;\tau)$, $\alpha \in [0,1), \beta \in [0,1], \tau \in \mathbb{R}^* = \mathbb{R}-\{0\}$. Then,

$$\sum_{n=2}^{\infty}|a_n| \leq \frac{(1-\alpha)|\tau|}{(1+\beta)\left[1+(1-\alpha)|\tau|\right]} \tag{2.7}$$

and

$$\sum_{n=2}^{\infty}n|a_n| \leq \frac{2(1-\alpha)|\tau|}{(1+\beta)\left[1+(1-\alpha)|\tau|\right]}. \tag{2.8}$$

**Poof.** Using Theorem 2.2, we obtain

$$\left[1+(1-\alpha)|\tau|\right](1+\beta)\sum_{n=2}^{\infty}|a_n| \leq \sum_{n=2}^{\infty}\left[n+(1-\alpha)|\tau|-1\right]\left[1+(n-1)\beta\right]|a_n| \leq (1-\alpha)|\tau|.$$

From here, we can easily show that (2.7) is true.
Similarly, we write

$$(1+\beta)\sum_{n=2}^{\infty}\left[n+(1-\alpha)|\tau|-1\right]|a_n| \leq \sum_{n=2}^{\infty}\left[n+(1-\alpha)|\tau|-1\right]\left[1+(n-1)\beta\right]|a_n| \leq (1-\alpha)|\tau|;$$

that is,

$$(1+\beta)\sum_{n=2}^{\infty}n|a_n| \leq (1-\alpha)|\tau|+\left[1-(1-\alpha)|\tau|\right](1+\beta)\sum_{n=2}^{\infty}|a_n|$$

Using (2.7) in the last inequality, we obtain

$$(1+\beta)\sum_{n=2}^{\infty}n|a_n| \leq \frac{2(1-\alpha)|\tau|}{1+(1-\alpha)|\tau|},$$

which immediately yields the inequality (2.8).

Thus, the proof of Theorem 2.3 is completed.

Setting $\tau = 1$ in Theorem 2.3, we obtain the following corollary.

**Corollary 2.7.** Let the function $f(z)$ definition by (1.2) belongs to the class $TS^*C(\alpha,\beta)$, $\alpha \in [0,1), \beta \in [0,1]$. Then,

$$\sum_{n=2}^{\infty}|a_n| \leq \frac{1-\alpha}{(2-\alpha)(1+\beta)}$$

and

$$\sum_{n=2}^{\infty}n|a_n| \leq \frac{2(1-\alpha)}{(2-\alpha)(1+\beta)}.$$

**Remark 2.5.** The result obtained in the Corollary 2.7 verifies to Lemma 2 (with $n = p = 1$) of [2].

From Theorem 2.2, for the coefficient bound estimates, we arrive at the following result.

**Corollary 2.8.** Let $f \in TS^*C(\alpha, \beta; \tau), \alpha \in [0,1), \beta \in [0,1], \tau \in \mathbb{R}^* = \mathbb{R} - \{0\}$. Then,

$$|a_n| \leq \frac{(1-\alpha)|\tau|}{[n+(1-\alpha)|\tau|-1][1+(n-1)\beta]}, \quad n = 2, 3, \ldots .$$

**Remark 2.6.** Numerous consequences of Corollary 2.8 can indeed be deduced by specializing the various parameters involved. Many of these consequences were proved by earlier workers on the subject (see, for example, [1, 9, 11]).

Acknowledgements

The author is grateful to the anonymous referees for their valuable comments and suggestions.